\documentclass[12pt]{article}
\usepackage{amssymb,latexsym,start2e,amsthm,amsmath}
\usepackage{pstricks,pst-node,pst-tree}
\usepackage{graphicx, pst-plot}

\newcommand{\Prob}{\operatorname{Prob}}

\newcommand{\Tha}{\hat{T}}
\newcommand{\n}[1]{\cnode*{.1}{#1}}
\newcommand{\mm}{\hspace*{-50pt}m{m\choose2}\hspace*{-50pt}}
\newcommand{\fone}{\frac{m-2}{3m}}
\newcommand{\ftwo}{\frac{m-1}{2m^2}}
\newcommand{\fthr}{\frac{m}{m^2}}
\newcommand{\ffou}{\frac{m}{m^3}}
\newcommand{\tone}{\frac{1}{2\cdot 3}}
\newcommand{\ttwo}{\frac{1}{2\cdot 1}}

\begin{document}
\pagestyle{plain}

\title{Probabilistic proofs of hook length formulas involving trees
}
\author{
Bruce E. Sagan\\[-5pt]
\small Department of Mathematics, Michigan State University,\\[-5pt]
\small East Lansing, MI 48824-1027, USA, \texttt{sagan@math.msu.edu}
}

\date{\today\\[10pt]
	\begin{flushleft}
	                                       \\[5pt]
{\it Dedication.\/}  This paper is dedicated to the memory of Pierre
Leroux.  In fact, reference~\cite{sy:pat} was written while both Sagan
and Yeh were visiting the Universit\'e de Qu\'ebex \`a Montr\'eal and
partaking of Leroux's legendary hospitality.
	\end{flushleft}
}

\maketitle

\begin{abstract}

Recently, Han discovered two formulas involving binary
trees which have the interesting property that hooklengths appear as
exponents.  The purpose of this note is to give a probabilistic proof
of one of Han's formulas.  Yang has generalized Han's results to
ordered trees.  We show how the probabilistic approach can also be
used in Yang's setting, as well as for a generalization of Han's
formula in terms of certain infinite trees.

\end{abstract}

\section{Introduction and definitions}

Frame, Robinson, and Thrall~\cite{frt:hgs} discovered the hook formula
for standard Young tableaux.  Greene, Nijenhuis, and
Wilf~\cite{gnw:ppf} then gave a probabilistic proof of this result
where the hook lengths appeared in a very natural way.  The same trio
also used probabilistic methods to prove the sum of squares formula
for the symmetric group~\cite{gnw:apm}.  Sagan~\cite{sag:srs} and
Sagan and Yeh~\cite{sy:pat} gave probabilistic algorithms for proving
hook formulas for shifted Young tableaux and trees, respectively.

Recently, inspired by an identity of Postnikov~\cite{pos:pab},
Han~\cite{han:nhl} proved two identities involving binary 
trees which have the interesting property that hooklengths appear as
exponents.  (Han~\cite{han:agp} also discovered an identity with this
same property which generalizes Postnikov's.)  Han's demonstration was
by algebraic manipulation of recursions.  Yang~\cite{yan:ghh}
generalized Han's identities to weighted ordered 
trees.  Again, the proofs were algebraic in nature, this time using
generating functions.  

The purpose of this note is to give a probabilistic proof of Han's
first formula which is similar in some ways to the second algorithm of
Greene, Nijenhuis, and Wilf.  
A weighted version of the algorithm 
proves the analogous identity of Yang.  A second generalization of
Han's original formula to certain infinite trees is also demonstrated by
this method.  The rest of this section is devoted to the necessary
definitions to state the identities to be proved.  Section~\ref{a}
gives the probabilistic algorithm and proofs.  The final section is
devoted to indicating how Han's second formula might be demonstrated
probabilistically.

\bfi
$$
\begin{psmatrix}[rowsep=15pt,colsep=25pt]
      &      &\n{a1}&&      &\n{b1}&&      &\n{c1}&      &&\n{d1}&      &&\n{e1}&      &\\
      &\n{a2}&      &&\n{b2}&      &&\n{c2}&      &\n{c3}&&      &\n{d2}&&      &\n{e2}&\\
\n{a3}&      &      &&      &\n{b3}&&      &      &      &&\n{d3}&      &&      &      &\n{e3}
\end{psmatrix}
\ncline{a1}{a3}
\ncline{b1}{b2}
\ncline{b2}{b3}
\ncline{c1}{c2}
\ncline{c1}{c3}
\ncline{d1}{d2}
\ncline{d2}{d3}
\ncline{e1}{e3}
$$
\vspace{5pt}
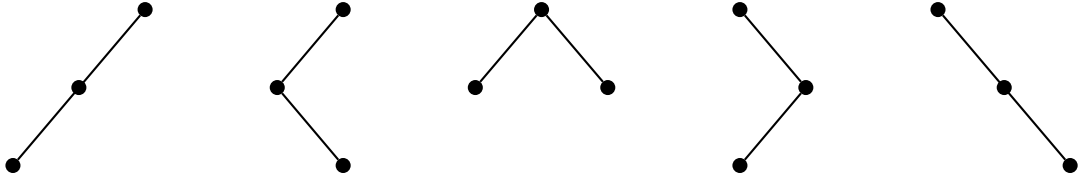
\capt{The trees in $\cB(3)$}
\label{B3}
\efi

For any tree, $T$, we denote the vertex set of $T$ by $V(T)$.
If no confusion will result, we
will often write $v\in T$ and $|T|$ in place of the more cumbersome
$v\in V(T)$ and $|V(T)|$, where $|\cdot|$ denotes cardinality.
If $T$ is rooted and $v\in T$,
then the set of children of $v$ will be denoted $C_v$, and we let
$c_v=|C_v|$.
The {\it hook  of $v$\/}, $H_v$, is the set of
descendents of $v$ (including $v$ itself) with corresponding {\it hook
length\/} $h_v=|H_v|$.

A {\it binary tree\/}, $T$, is a rooted tree where every vertex has
either no children, a left child, a right child, or both children.
Let $\cB$ denote the family of all binary trees and 
let 
$$
\cB(n)=\{T\in \cB\ :\ |T|=n\}.
$$
For example, the trees in $\cB(3)$ are displayed in Figure~\ref{B3}.
In what follows, we will use similar notation for other families of
trees.  The formula of Han which we will prove is as follows.
\bth[Han]
For each positive integer $n$ we have
\beq
\label{H}
\sum_{T\in \cB(n)} \prod_{v\in T} \frac{1}{h_v 2^{h_v-1}} =\frac{1}{n!}.
\eeq
\eth

\bfi
$$
\begin{psmatrix}[rowsep=15pt,colsep=23pt]
\n{a1}&&      &\n{b1}&      &&      &\n{c1}&      &&      &\n{d1}&      &&      &\n{e1}&\\
\n{a2}&&      &\n{b2}&      &&\n{c2}&      &\n{c3}&&\n{d2}&      &\n{d3}&&\n{e2}&\n{e3}&\n{e4}\\
\n{a3}&&\n{b3}&      &\n{b4}&&\n{c4}&      &      &&      &      &\n{d4}&&      &\\
\n{a4}\\
m^3   &&      &\mm   &      &&      &\mm   &      &&      &\mm   &      &&      &{m\choose3}
\end{psmatrix}
\ncline{a1}{a4}
\ncline{b1}{b2}
\ncline{b2}{b3}
\ncline{b2}{b4}
\ncline{c1}{c2}
\ncline{c1}{c3}
\ncline{c2}{c4}
\ncline{d1}{d2}
\ncline{d1}{d3}
\ncline{d3}{d4}
\ncline{e1}{e2}
\ncline{e1}{e3}
\ncline{e1}{e4}
$$
\vspace{5pt}
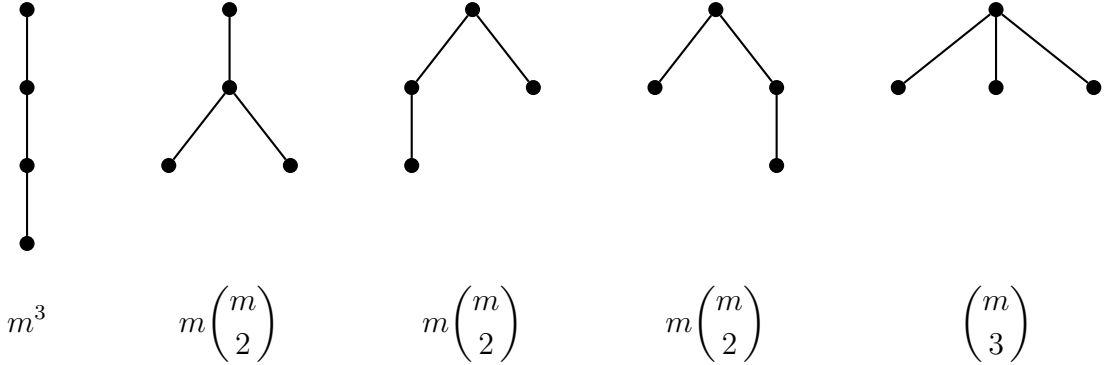
\capt{The trees in $\cO(4)$ and their weights}
\label{O4}
\efi

Now consider finite ordered trees weighted by
$$
w(T)=\prod_{v\in T} {m\choose c_v}
$$
where $m$ is a variable.  Let $\cO$ denote the family of weighted
ordered trees.  
The trees in $\cO(4)$ along with their weights are shown in Figure~\ref{O4}.
Then the identity of Yang we are considering is
equivalent to:
\bth[Yang]
For each positive integer $n$ we have
\beq
\label{Y}
\sum_{T\in \cO(n)}w(T) \prod_{v\in T} \frac{1}{h_v m^{h_v-1}} =\frac{1}{n!}.
\eeq
\eth

Some comments about this result are in order.  First of all, it is
remarkable because the right-hand side of the equation does not depend
on $m$.  Secondly, it implies Han's formula by letting $m=2$, since
then $w(T)$ is just the number of ways of turning an ordered tree into
a binary tree.  Finally, Yang's weighting was actually
$$
w(T)=\prod_{v\in T} {m\choose c_v} s^{c_v}= s^{|T|-1}\prod_{v\in T} {m\choose c_v}
$$
where $s$ is another parameter.  So one can recover Yang's equation by
multiplying both sides of~\ree{Y} by $s^{n-1}$.  Also, Yang assumes
that $m$ and $s$ are constants satisfying certain conditions, but it
is clearly not necessary to do so.

\bfi
$$
\begin{psmatrix}[rowsep=15pt,colsep=25pt]
      &      &      &&      &      &       &      &\n{t1}&      &      &      &      &&      &      \\
      &      &      &&      &\Tb=  &       &\n{t2}&      &\n{t3}&      &      &      &&      &      \\
      &      &      &&      &      &\n{t4} &\n{t5}&\n{t6}&\n{t7}&      &      &      &&      &      \\
      &      &      &&      &      &\vdots &\vdots&\vdots&\vdots&      &      &      &&      &      \\[20pt]
      &      &\n{a1}&&      &\n{b1}&       &      &\n{c1}&      &      &\n{d1}&      &&\n{e1}&      \\
      &\n{a2}&      &&\n{b2}&      &       &\n{c2}&      &      &\n{d2}&      &\n{d3}&&      &\n{e2}\\
\n{a3}&      &      &&\n{b3}&      &       &      &\n{c3}&      &      &      &      &&      &\n{e3}
\end{psmatrix}
\ncline{a1}{a3}
\ncline{b1}{b2}
\ncline{b2}{b3}
\ncline{c1}{c2}
\ncline{c2}{c3}
\ncline{d1}{d2}
\ncline{d1}{d3}
\ncline{e1}{e2}
\ncline{e2}{e3}
\ncline{t1}{t2}
\ncline{t1}{t3}
\ncline{t2}{t4}
\ncline{t2}{t5}
\ncline{t2}{t6}
\ncline{t3}{t7}
$$
\vspace{5pt}
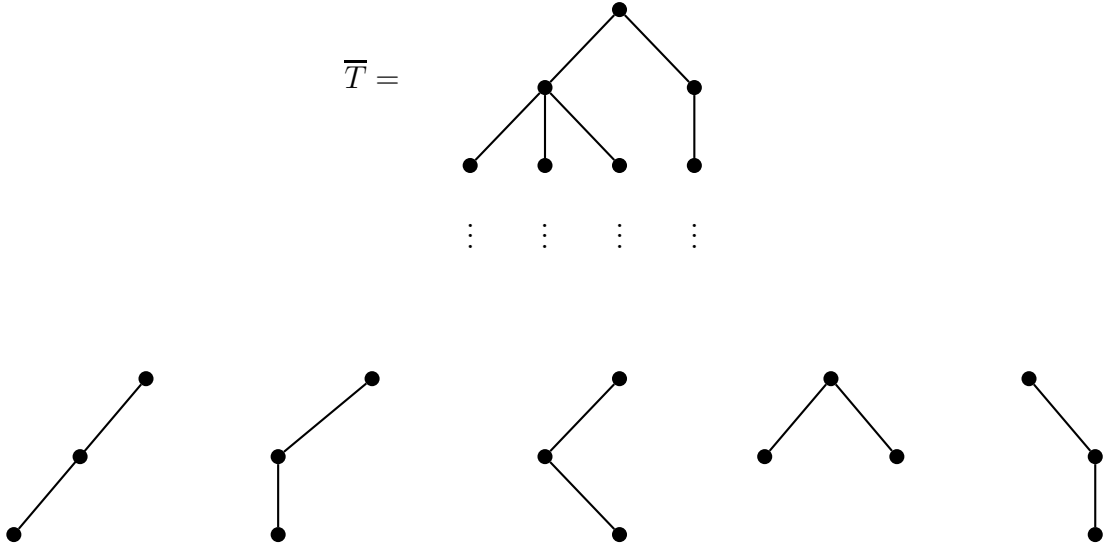
\capt{The subtrees in $\cTb(3)$ for the given $\Tb$}
\label{T3}
\efi

For our second generalization of equation~\ree{H}, consider a fixed,
infinite, ordered tree $\Tb$ such that $0<\cb_v<\infty$ for all $v\in
T$.  We are using $\cb_v$ for the number of children of $v$ to
emphasize that this is being calculated in $\Tb$.
Let $\cTb$ be the family of all subtrees of $\Tb$.  Since $\Tb$ is
ordered, its vertices are distinguishable, i.e., $V(\Tb)$ is a set
rather than a multiset.  So we consider two subtrees $T,T'$ to be
equal if and only if $V(T)=V(T')$.  For example, $\cTb=\cB$ if we
let $\Tb$ be the tree with $\cb_v=2$ for all $v\in \Tb$.  As another
illustration, Figure~\ref{T3} shows part of one possible $\Tb$ and the
corresponding trees in $\cTb(3)$.
\bth
For each positive integer $n$ and each tree $\Tb$ satisfying the above
restrictions, we have
\beq
\label{S}
\sum_{T\in \cTb(n)} \prod_{v\in T} \frac{1}{h_v \cb_v^{h_v-1}}=\frac{1}{n!}.
\eeq
\eth

\section{The algorithm}
\label{a}

For any rooted tree $T$, an {\it increasing labeling of $T$\/} is a
bijection 
$$
\ell:T\ra\{1,2,\ldots,|T|\}
$$ 
such that for any $v\in T$ and
any $w\in C_v$ we have $\ell(v)<\ell(w)$.  We will often let $L=L(T)$
stand for an increasing labeling of $T$ viewed as $T$ with the labels
attached to its vertices.  Similarly, we will write $L(\cF)$ for the
family of all increasing labelings of trees in the family $\cF$.
Let $f^T$ be the number of
increasing labelings $L(T)$ where $T$ has distinguishable vertices.
The following hook length formula for $f^T$ is well known and easy to
prove
\beq
\label{fT}
f^T=\frac{n!}{\prod_{v\in T} h_v}.
\eeq
So if we multiply any of the three identities from the previous
section by $n!$, we obtain a sum of the form
$$
\sum_{T\in\cF(n)} f^T \pi(T) = 1
$$
where $\cF$ is a family of trees and $\pi(T)$ is a product.
We wish to interpret $\pi(T)$ as the
probability of obtaining an increasing labeling $L$ of $T$ for an
appropriate probability distribution on $L(\cF(n))$.  The identity
will then follow since
$$
1=\sum_{L\in L(\cF(n))} \Prob(L) 
= \sum_{T\in\cF(n)} \sum_{L=L(T)} \pi(T)
= \sum_{T\in\cF(n)} f^T \pi(T).
$$
Note that $\Prob(L)$ will depend on $T$ where $L=L(T)$, but not on the
specific labeling of $T$.

The probability distribution will be obtained by
specifying the parameters in the following basic algorithm which takes
as input a positive integer $n$ and a family of trees $\cF$ and
outputs a labeling $L$ of some $T\in\cF(n)$.  
\ben
\item Let $L$ consist of a single root labeled $1$.
\item While $|L|<n$, consider all possible leaves $v$ one could add to
  $L$ and still stay in $L(\cF)$.  Pick one such leaf, label
  it $|L|+1$, and add it to $L$ with probability $\Prob(v,L)$. 
\item Output $L$.
\een
It will be convenient to also use the notation $\Prob(v,L)$ when 
$v\in L$.  In that case, it should be interpreted as $\Prob(v,L')$
where $L'$ is subtree of $L$ induced by those vertices with labels
less than $\ell(v)$.

To finish the proofs, we just need to specify for each of the three
families what the probabilities $\Prob(v,L)$ are, and prove that they
describe a probality distribution such that all increasing labelings $L$
of a given tree are equally likely with the common value being 
$$
\Prob(L)=\prod_{v\in L} \Prob(v,L) =\pi(T).
$$

\medskip

{\it Proof of~\ree{H}.\/}
Given a tree $T$ rooted at $r$ and $v\in T$ we let $P_v$ be the unique
$r\hor v$ path.  The {\it depth of $v$\/}, $d_v$, is the length of $P_v$. 
In the algorithm, let
$$
\Prob(v,L)=\frac{1}{2^{d_v}}.
$$
For example, Figure~\ref{ProbB} shows one of the trees $T$ in $\cB(3)$
along with these probabilities for each possible leaf $v$ which could
be added to $T$.  To further distinguish such leaves from the nodes of $T$,
the corresponding edges are dashed.

\bfi
$$
\begin{psmatrix}[rowsep=10pt,colsep=30pt]
      &       & \n{a} &       \\
      &       &       &       \\
      & \n{b} &       & \n{c} \\
      &\st{20}&       & 1/2   \\
\n{d} &       & \n{e} &       \\
 1/4  &       &       &       \\
      & \n{f} &       & \n{g} \\
      & 1/8   &       & 1/8      
\end{psmatrix}
\ncline{a}{b}
\ncline[linestyle=dashed]{a}{c}
\ncline[linestyle=dashed]{b}{d}
\ncline{b}{e}
\ncline[linestyle=dashed]{e}{f}
\ncline[linestyle=dashed]{e}{g}
$$
\vspace{5pt}
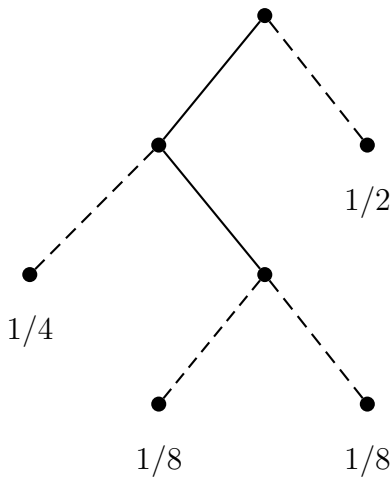
\capt{A tree in  $\cB(3)$ and the probabilities of its additional leaves}
\label{ProbB}
\efi

We first need a lemma which will be used in all three proofs.  So we
will state it in general terms.
\ble
\label{F}
For each of the three families $\cF$ under consideration and
for each $L\in L(\cF)$ we have
\beq
\label{sumv}
\sum_{v} \Prob(v,L)=1
\eeq
where the sum is over all leaves $v$ which could be added to $L$.
\ele
{\it Proof of the Lemma for $\cF=\cB$.\/}  We induct on $|L|$ where
the base case is easy to do.  Given $L$, let $w$ be the leaf of $L$
such that $\ell(w)=|L|$ and let $L'=L-w$.
Then the
terms in the sum for $L'$ are the same as those in the sum for $L$
except that the summand $1/2^{d_w}$ in the former has been replaced by 
$1/2^{d_w+1}+1/2^{d_w+1}$.  Since these two expressions are equal, so
are the sums, and we are done by induction.\hqedm

Next we need to verify that for $L=L(T)$ we have $\Prob(L)=\pi(T)$, i.e.,
$$
\Prob(L)=\prod_{v\in T}\frac{1}{2^{h_v-1}}
$$
Again, let $\ell(w)=|L|$ and $L'=L-w$.  Then the hook lengths in $L$
are the same as those in $L'$ except that the $d_w$ vertices on
the path $P_w-w$ 
have all been increased by one.  Note also that $w$ itself does not
contribute to the product above since $1/2^{h_w-1}=1$.  So, by induction,
\beq
\label{B}
\Prob(L) =\Prob(w,L')\Prob(L')
=\frac{1}{2^{d_w}}\prod_{v'\in T'}\frac{1}{2^{h_{v'}-1}}
=\prod_{v\in T}\frac{1}{2^{h_v-1}}.
\eeq

There remains to show that $\Prob(L)$ defines a probability
distribution.  But using the Lemma and induction as well as our usual
notation:
\bea
\sum_{L\in L(\cB(n))} \Prob(L)
&=&\dil\sum_{L\in L(\cB(n))} \Prob(w,L')\Prob(L')\\[10pt]
&=&\dil\sum_{L'\in L(\cB(n-1))} \Prob(L') \sum_w \Prob(w,L')\\[10pt]
&=&\dil\sum_{L'\in L(\cB(n-1))} \Prob(L')\\[10pt]
&=&1.
\eea
This finishes the proof of~\ree{H}.\hqedm

Note that the proof that $\Prob(L)$ forms a probability distribution
only depends on Lemma~\ref{F}.  So in the next two proofs, we will
skip this step.

\medskip

{\it Proof of~\ree{Y}.}   Note that the left-hand side of~\ree{Y} is a
rational function of $m$, so clearing denominators gives a polynomial
equation.  Thus it suffices to prove that this identity holds for
infinitely many values of $m$.  We will provide a proof for all real numbers
$m\ge M$ where $M$ is chosen sufficiently large such that 
$0\le \Prob(v,L)\le 1$ for all $v\in L\in L(\cO(n))$.  This will be
possible because $\Prob(v,L)$ will be asymptotic to, but smaller than
or equal to, $1/m^{d_v-1}$.   Specifically, let
$$
\Prob(v,L)=\frac{m-c_p}{(c_p+1)m^{d_v}}
$$
where $p$ is the parent of $v$.  Remember that, according to our
convention following the description of the algorithm, $c_p$ is
calculated in the subtree of $L$ induced by those vertices with labels
less that $\ell(v)$.  In particular, $c_p$ does not count $v$ itself.
Figure~\ref{ProbO} displays a tree of $\cO(4)$ and the probabilities
of the leaves which can be added to it.

\bfi
$$
\begin{psmatrix}[rowsep=10pt,colsep=50pt]
      &       & \n{a} &       &        \\
      &       &       &       &        \\
      &       &       &       &        \\
\n{b} & \n{c} & \n{d} & \n{e} & \n{f}  \\
\fone &       & \fone &       & \fone  \\
      &       &       &       &        \\
      & \n{g} & \n{h} & \n{i} & \n{j}  \\
      & \fthr & \ftwo &       & \ftwo  \\
      &       &       &       &        \\
      &       &       & \n{k} &        \\
      &       &       & \ffou &      
\end{psmatrix}
\ncline[linestyle=dashed]{a}{b}
\ncline{a}{c}
\ncline[linestyle=dashed]{a}{d}
\ncline{a}{e}
\ncline[linestyle=dashed]{a}{f}
\ncline[linestyle=dashed]{c}{g}
\ncline[linestyle=dashed]{e}{h}
\ncline{e}{i}
\ncline[linestyle=dashed]{e}{j}
\ncline[linestyle=dashed]{i}{k}
$$
\vspace{5pt}
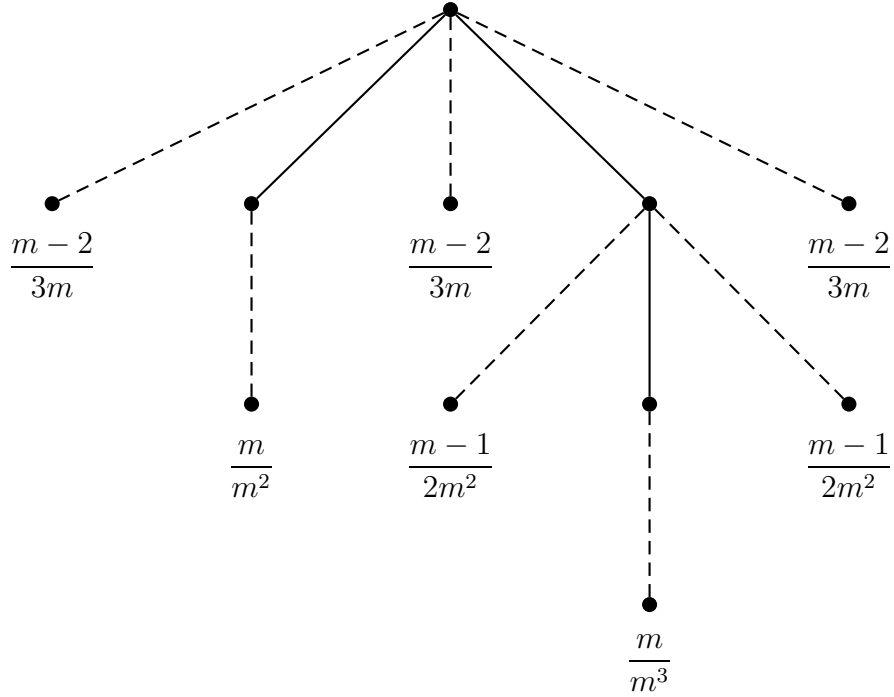
\capt{A tree in  $\cO(4)$ and the probabilities of its additional leaves}
\label{ProbO}
\efi

Our first order of business will be to prove Lemma~\ref{F} in this
setting.

{\it Proof of the Lemma for $\cF=\cO$.\/}  As before, we induct on
$L$, keeping the notation the same as the first proof.  We also let
$p$ be the parent of $w$ and write $p'$ if we are considering $p$ as a
vertex of $L'$. So $c_p=c_{p'}+1$ and the terms in the sum for $L'$
corresponding to the $c_{p'}+1$ possible children which could be added to
$p'$ give a total of
$$
(c_{p'}+1)\frac{m-c_{p'}}{(c_{p'}+1)m^{d_w}}=\frac{m-c_p+1}{m^{d_w}}.
$$
In the sum for $L$, these terms are replaced by $c_p+1$ summands for
children of $p$ and one for a child of $w$, giving
$$
(c_p+1)\frac{m-c_p}{(c_p+1)m^{d_w}}+\frac{m}{m^{d_w+1}}=\frac{m-c_p+1}{m^{d_w}}.
$$
Since these two contributions are the same and all other terms in two
sums match up, we are done.\hqedm

We next need to show that $\Prob(L)=\pi(T)$ for $\cF=\cO$.  Keeping
our usual notation we have $\Prob(L)/\Prob(L')=\Prob(w,L')$.  So the
desired equality will follow by induction, the reasoning applied to
obtain~\ree{B}, and the computation
$$
\frac{\pi(T)}{\pi(T')}=
\frac{\prod_{v\in T}{m\choose c_v}/m^{h_v-1}}{\prod_{v'\in T'}{m\choose c_{v'}}/m^{h_{v'}-1}}
=\frac{{m\choose c_p}{m\choose c_w}}{{m\choose c_{p'}} m^{d_w}}
=\frac{m-c_{p'}}{(c_{p'}+1)m^{d_w}}=\Prob(v,L').
$$
This completes the proof for $\cO$.\hqedm

\bfi
$$
\begin{psmatrix}[rowsep=10pt,colsep=30pt]
       &      &\n{t1}&      \\
       &      &      &      \\
       &\n{t2}&      &\n{t3}\\
       &      &      &      \\
\n{t4} &\n{t5}&\n{t6}&\n{t7}\\
\tone  &\tone &\tone &\ttwo
\end{psmatrix}
\ncline{t1}{t2}
\ncline{t1}{t3}
\ncline[linestyle=dashed]{t2}{t4}
\ncline[linestyle=dashed]{t2}{t5}
\ncline[linestyle=dashed]{t2}{t6}
\ncline[linestyle=dashed]{t3}{t7}
$$
\vspace{5pt}
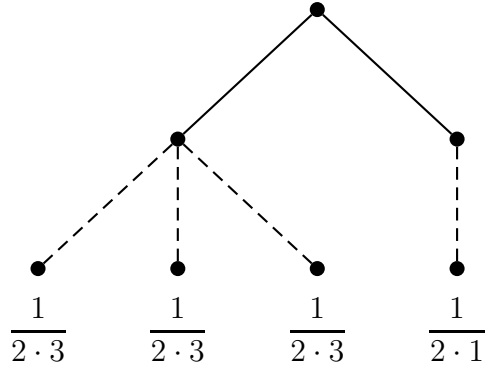
\capt{A subtree in $\cTb(3)$ and additional leaves}
\label{ProbT}
\efi

{\it Proof of~\ree{S}.\/}  
For this case, we proceed as usual, but letting
$$
\Prob(v,L)=\prod_{x\in P_v-v} \frac{1}{\cb_x}.
$$
Figure~\ref{ProbT} gives an example using a tree from $\cTb(3)$ where
$\cTb$ is as in Figure~\ref{T3}.

{\it Proof of the Lemma for $\cF=\cTb$.\/} Now in passing from the sum
for $L'$ to the sum for $L$, a single term 
$\prod_{x\in P_w-w} 1/\cb_x$ has been replaced by $\cb_w$ terms all
equal to $\prod_{x\in P_w} 1/\cb_x$.  Clearly this does not change the
sum.\hqedm

The proof that $\Prob(L)=\pi(t)$ is just like the one for $\cB$ except
that the hook length powers of 2 are replaced by powers of $\cb_x$.
So we are done with the case of $\cTb$.\hqedm

\section{An open problem}

As remarked in the introduction, Han actually proved two formulas
in~\cite{han:nhl}, both having hook lengths as exponents.  We have
unable to give a probabilistic proof of the second one.  But will
indicate how one might go in the hopes that someone else may be
able to push it through.

Call a binary tree {\it complete\/} if every vertex has 0 or 2
children.  Given a binary tree $T$ on $n$ nodes it has 
{\it completion $\Tha$\/} which is the complete binary tree obtained from $T$
by adding all $n+1$ possible leaves.  If $T$ is the tree with the
solid edges in Figure~\ref{ProbB} then $\Tha$ is the tree which also
includes the dashed edges.
It is not hard to show
using~\ree{fT} that
\beq
\label{fTh}
f^{\Tha}=\frac{(2n+1)!}{\prod_{v\in T} (2h_v+1)}.
\eeq

Han's second formula is
$$
\sum_{T\in\cB(n)}\prod_{v\in T} \frac{1}{(2h_v+1)2^{2h_v-1}}=\frac{1}{(2n+1)!}.
$$
Using~\ree{fTh}, this can be rewritten as
$$
\sum_{T\in\cB(n)}f^{\Tha}\prod_{v\in T} \frac{1}{2^{2h_v-1}}=1.
$$
It would be very interesting to find a probability distribution on increasing labelings of
complete trees $\Tha$ whose probabilities are given by $\prod_{v\in T}
1/2^{2h_v-1}$.  Once this is done, similar ideas should prove the
generalization to $\cO$ due to Yang~\cite{yan:ghh}.  It is not clear
how to generalize Han's formula to the $\cTb$ case, but would be
interesting to do if possible.

\bigskip
\bibliographystyle{acm}
\begin{small}
\bibliography{ref}
\end{small}

\end{document}